\documentclass[10pt,a4paper]{article}
\usepackage{amsfonts}
\usepackage{amsmath}
 \title{\textbf{$3$-Sasakian manifolds, $3$-cosymplectic manifolds and Darboux theorem}}
 \author{\begin{tabular}{cc}
  \textsc{Beniamino Cappelletti Montano,  Antonio De Nicola}\\
  Department of Mathematics,  University of Bari \\
  Via E. Orabona 4,\\
  I-70125 Bari (Italy) \\
  {\texttt{cappelletti@dm.uniba.it,}\  \texttt{antondenicola@gmail.com}}
   \end{tabular}}
   \date{}
\usepackage{graphicx}
\usepackage{amsmath}
\newtheorem{theorem}{Theorem}[section]

\newtheorem{corollary}[theorem]{Corollary}

\newtheorem{lemma}[theorem]{Lemma}

\newtheorem{proposition}[theorem]{Proposition}
\newtheorem{remark}[theorem]{Remark}

\newenvironment{proof}[1][Proof]{\textbf{#1.} }{\ \rule{0.5em}{0.5em}}

\begin{document}

\maketitle

\begin{center}
\textit{Dedicated to the memory of Giulio Minervini on the
anniversary of his departure.}
\end{center}

\begin{abstract}
We present a compared analysis  of some properties of $3$-Sasakian
and $3$-cosymplectic manifolds. We construct a canonical connection on
an almost $3$-contact metric manifold which generalises the Tanaka-Webster
connection of a contact metric manifold and we use this connection to show that
a $3$-Sasakian manifold does not admit any Darboux-like coordinate system.
Moreover, we prove that any $3$-cosymplectic manifold is Ricci-flat and admits
a Darboux coordinate system if and only it is flat.
\end{abstract}

\section{Introduction}
Both $3$-Sasakian and $3$-cosymplectic manifolds belong to the class of
almost contact (metric) $3$-structures, introduced by Kuo (\cite{kuo}) and,
independently, by Udriste (\cite{udriste}).
The study of $3$-Sasakian manifolds has been conducted by several authors (see for example \cite{galicki1,galicki2}
and references therein) due to the increasing awareness of their importance in mathematics and in physics,
together with the closely linked hyper-K\"{a}hlerian and quaternionic K\"{a}hlerian manifolds.
Recently they have made an appearance also in supergravity and M-theory (see \cite{acharya,agricola,gibbons}).
Less studied, so far, are $3$-cosymplectic manifolds, also called hyper-cosymplectic, but we can
list some recent publications (\cite{pastore, kashiwada2, martin,tshikuna-matamba}). For example, Kashiwada and his collaborators proved in \cite{kashiwada2} that any $b$-Kenmotsu (see \cite{blair1,janssens}) almost contact $3$-structure must be $3$-cosymplectic.

In this paper we present a compared analysis  of some properties of $3$-Sasakian
and $3$-cosymplectic manifolds. We start with a brief review of some known results
on these classes of manifolds, contained in Section 2. In the third section we construct
a canonical connection on an almost $3$-contact metric manifold and we study its curvature and torsion
analysing also its behaviour in the special cases of $3$-Sasakian and $3$-cosymplectic manifolds.
Our connection can be interpreted as a generalisation of the (generalised) Tanaka-Webster connection of a contact metric manifold, introduced by Tanno in \cite{tanno}. The section is concluded by a further investigation of the properties of $3$-cosymplectic manifolds concerning with their projectability which leads us to prove that every $3$-cosymplectic manifold is Ricci-flat.
In the final section we analyse the possibility of finding a Darboux-like coordinate system
on $3$-Sasakian and $3$-cosymplectic manifolds. Firstly we establish a relation which
holds in any almost $3$-contact metric manifold linking the horizontal part of the metric with the three fundamental forms $\Phi_\alpha$. This relation is responsible for a kind of rigidity of this class of manifolds which links
the existence of Darboux coordinates to the flatness of the manifold and does not hold in the case of a single Sasakian or cosymplectic structure. In particular, on the one hand, using our canonical connection and the (local) projection of  a $3$-Sasakian manifold over a quaternionic K\"{a}hlerian manifold (see \cite{galicki1,ishihara}), we show that $3$-Sasakian manifolds, unlike the Sasakian ones, do not admit any Darboux-like coordinate system. This result is related to the fact that $3$-Sasakian manifolds are not (horizontally) flat. On the other hand, we show that a $3$-cosymplectic manifold admits a Darboux coordinate system in the neighbourhood of each point if and only its metric is flat.
\\
\\
\textbf{Acknowledgements} \ The authors wish to thank Prof. A. M. Pastore for the helpful suggestions during the preparation of the paper.

\section{Preliminaries}
An \emph{almost contact manifold} is an odd-dimensional manifold
$M$ which carries a field $\phi$ of endomorphisms of the tangent
spaces, a vector field $\xi$, called \emph{characteristic} or
\emph{Reeb vector field}, and a $1$-form $\eta$ satisfying
$\phi^2=-I+\eta\otimes\xi$ and  $\eta\left(\xi\right)=1$, where $I
\colon TM\rightarrow TM$ is the identity mapping. From the
definition it follows also that $\phi\xi=0$, $\eta\circ\phi=0$ and
that the $(1,1)$-tensor field $\phi$ has constant rank $2n$ (cf.
\cite{blair1}). An almost contact  manifold
$\left(M,\phi,\xi,\eta\right)$ is said to be \emph{normal} when
the tensor field $N=\left[\phi,\phi\right]+2d\eta\otimes\xi$
vanishes identically, $\left[\phi,\phi\right]$ denoting the
Nijenhuis tensor of $\phi$. It is known that any almost contact
manifold $\left(M,\phi,\xi,\eta\right)$ admits a Riemannian metric
$g$ such that
\begin{equation}\label{compatibile}
g\left(\phi E,\phi F\right)=g\left(E,F\right)-\eta\left(E\right)\eta\left(F\right)
\end{equation}
holds for all $E,F\in\Gamma\left(TM\right)$. This metric $g$ is called a
\emph{compatible metric} and the manifold $M$ together with the structure
$\left(\phi,\xi,\eta,g\right)$ is called an \emph{almost contact
metric manifold}. As an immediate consequence of
\eqref{compatibile}, one has
$\eta=g\left(\cdot,\xi\right)$. The $2$-form $\Phi$ on $M$
defined by $\Phi\left(E,F\right)=g\left(E,\phi F\right)$ is called
the \emph{fundamental $2$-form} of the almost contact metric manifold
$M$.  Almost contact metric manifolds such that
both $\eta$ and $\Phi$ are closed are called \emph{almost
cosymplectic manifolds} and almost contact metric manifolds such that $d\eta=\Phi$ are called \emph{contact metric manifolds}. Finally, a normal almost cosymplectic manifold
is called a \emph{cosymplectic manifold} and a normal contact metric manifold is said to be a \emph{Sasakian manifold}. In terms of the covariant derivative of $\phi$ the cosymplectic and the Sasakian conditions can be expressed respectively by
\begin{equation*}
\nabla\phi=0
\end{equation*}
and
\begin{equation*}
\left(\nabla_E\phi\right)F=g\left(E,F\right)\xi-\eta\left(F\right)E
\end{equation*}
for all $E,F\in\Gamma\left(TM\right)$. It should be noted that both in Sasakian and in cosymplectic manifolds $\xi$ is a Killing vector field.

An  \emph{almost $3$-contact manifold}  is a
$\left(4n+3\right)$-dimensional smooth  manifold $M$
endowed with three almost contact structures
$\left(\phi_1,\xi_1,\eta_1\right)$,
$\left(\phi_2,\xi_2,\eta_2\right)$,
$\left(\phi_3,\xi_3,\eta_3\right)$ satisfying the following
relations, for every $\alpha,\beta\in\left\{1,2,3\right\}$,
\begin{equation}
\phi_\alpha\phi_\beta-\eta_\beta\otimes\xi_\alpha = \sum_{\gamma=1}^{3}\epsilon_{\alpha\beta\gamma}\phi_\gamma - \delta_{\alpha\beta}I,\quad
\phi_\alpha \xi_\beta = \sum_{\gamma=1}^{3}\epsilon_{\alpha\beta\gamma} \xi_\gamma,\quad \eta_\alpha\circ\phi_\beta =  \sum_{\gamma=1}^{3}\epsilon_{\alpha\beta\gamma} \eta_\gamma, \label{3-sasaki}
\end{equation}
where $\epsilon_{\alpha\beta\gamma}$ is the totally antisymmetric
symbol. This notion was introduced  by Kuo (\cite{kuo}) and,
independently, by Udriste (\cite{udriste}). In \cite{kuo} Kuo
proved that given an almost contact $3$-structure
$\left(\phi_\alpha,\xi_\alpha,\eta_\alpha\right)$, there exists a
Riemannian metric $g$ compatible with each of them and hence we
can speak of \emph{almost contact metric $3$-structures}. It is
well known that in any almost $3$-contact metric manifold the Reeb
vector fields $\xi_1,\xi_2,\xi_3$ are orthonormal with respect to
the compatible metric $g$ and that the structural group of the
tangent bundle is reducible to $Sp\left(n\right)\times I_3$.
Moreover, by putting
${\cal{H}}=\bigcap_{\alpha=1}^{3}\ker\left(\eta_\alpha\right)$ one
obtains a $4n$-dimensional distribution on $M$ and the tangent
bundle splits as the orthogonal sum
$TM={\cal{H}}\oplus\left\langle\xi_1,\xi_2,\xi_3\right\rangle$.
For a reason which will be clearer later we call any vector
belonging to the distribution $\cal H$ "horizontal" and any vector
belonging to the distribution
$\left\langle\xi_1,\xi_2,\xi_3\right\rangle$ "vertical". An almost
$3$-contact manifold $M$  is  said  to  be \emph{hyper-normal}  if
each  almost contact  structure
$\left(\phi_\alpha,\xi_\alpha,\eta_\alpha\right)$ is normal.

When the three structures
$\left(\phi_\alpha,\xi_\alpha,\eta_\alpha,g\right)$ are contact
metric structures, we say that $M$ is a \emph{$3$-contact
metric manifold} and when they are Sasakian, that is when each
structure $\left(\phi_\alpha,\xi_\alpha,\eta_\alpha\right)$ is
also normal, we call $M$ a \emph{$3$-Sasakian manifold}.
However these two notions coincide. Indeed as it has been proved
in 2001 by Kashiwada (\cite{kashiwada}),  every contact metric
$3$-structure is $3$-Sasakian. In any $3$-Sasakian manifold we have that, for each
$\alpha\in\left\{1,2,3\right\}$,
\begin{equation}\label{sasaki}
    \phi_\alpha=-\nabla\xi_\alpha.
\end{equation}
Using this, one obtains that $\left[\xi_1,\xi_2\right]=2\xi_3$,
$\left[\xi_2,\xi_3\right]=2\xi_1$,
$\left[\xi_3,\xi_1\right]=2\xi_2$. In particular, the vertical
distribution $\left\langle\xi_1,\xi_2,\xi_3\right\rangle$ is
integrable and defines a $3$-dimensional foliation of $M$
denoted by ${\cal{F}}_{3}$. Since $\xi_1$, $\xi_2$, $\xi_3$ are
Killing vector fields, ${\cal{F}}_{3}$ is a Riemannian foliation.
Moreover it has totally geodesic leaves of constant curvature $1$.
On the contrary, in a $3$-Sasakian manifold the horizontal
distribution ${\cal{H}}$ is never integrable.
About the foliation ${\cal F}_3$, Ishihara (\cite{ishihara}) has shown that
if ${\cal F}_3$ is regular then the space of leaves is a quaternionic
K\"{a}hlerian manifold.  Boyer, Galicki and  Mann have proved
the following more general result.

\begin{theorem}[\cite{galicki1}]\label{proiezione}
Let $\left(M^{4n+3},\phi_\alpha,\xi_\alpha,\eta_\alpha,g\right)$
be a $3$-Sasakian manifold such that the Killing vector fields
$\xi_1$, $\xi_2$, $\xi_3$ are complete. Then
\begin{description}
    \item[(i)]
    $M^{4n+3}$ is
    an Einstein manifold of positive scalar curvature equal to
    $2\left(2n+1\right)\left(4n+3\right)$.
    \item[(ii)] Each leaf $\cal L$ of the foliation ${\cal F}_3$
    is a $3$-dimensional homogeneous spherical space form.
    \item[(iii)] The space of leaves $M^{4n+3}/\cal F$ is a quaternionic
    K\"{a}hlerian orbifold of dimension $4n$ with positive scalar
    curvature equal to $16n\left(n+2\right)$.
\end{description}
\end{theorem}

By  an   \emph{almost $3$-cosymplectic manifold}   we  mean an
almost  $3$-contact  metric  manifold $M$ such that each almost
contact metric structure
$\left(\phi_\alpha,\xi_\alpha,\eta_\alpha,g\right)$ is almost
cosymplectic. The almost $3$-cosymplectic structure
$\left(\phi_\alpha,\xi_\alpha,\eta_\alpha,g\right)$ is called
\emph{$3$-cosymplectic} if it is hyper-normal. In this case $M$ is
said to be a \emph{$3$-cosymplectic manifold}. However it has been
proved recently that these two notions are the same:

\begin{theorem}\textbf{\emph{(\cite[Theorem 4.13]{pastore})}}
Any almost $3$-cosymplectic manifold is $3$-cosymplectic.
\end{theorem}

In any $3$-cosymplectic manifold we have that $\xi_\alpha$,
$\eta_\alpha$, $\phi_\alpha$ and $\Phi_\alpha$ are
$\nabla$-parallel. In particular
\begin{equation}\label{commutatore}
    \left[\xi_\alpha,\xi_\beta\right] = \nabla_{\xi_\alpha}\xi_\beta-\nabla_{\xi_\beta}\xi_\alpha = 0
\end{equation}
for all $\alpha,\beta\in\left\{1,2,3\right\}$, so that, as in any
$3$-Sasakian manifold,
$\left\langle\xi_1,\xi_2,\xi_3\right\rangle$ defines a
$3$-dimensional foliation ${\cal F}_3$ of $M^{4n+3}$.
However, unlike the case of $3$-Sasakian geometry, the horizontal subbundle
$\cal H$ of a $3$-cosymplectic manifold is integrable because,
for all $X,Y\in\Gamma\left(\cal H\right)$,
$\eta_\alpha\left(\left[X,Y\right]\right)=-2d\eta_\alpha\left(X,Y\right)=0$ since $d\eta_\alpha=0$.

\section{Further properties of $3$-Sasakian and $3$-cosymplectic manifolds}

In this section we investigate on further properties of $3$-Sasakian and $3$-cosymplectic manifolds. We start with the following preliminary result.

\begin{lemma}
Let $\left(M,\phi_\alpha,\xi_\alpha,\eta_\alpha,g\right)$
be an almost $3$-contact metric manifold. Then if $M$ is
$3$-Sasakian we have, for each $\alpha,\beta\in\left\{1,2,3\right\}$,
\begin{equation}\label{derivatalie1}
    {\cal L}_{\xi_\alpha}\phi_\beta=2\sum_{\gamma=1}^{3}\epsilon_{\alpha\beta\gamma}\phi_\gamma,
\end{equation}
and if $M$ is $3$-cosymplectic,
\begin{equation}\label{derivatalie2}
{\cal L}_{\xi_\alpha}\phi_\beta=0.
\end{equation}
\end{lemma}
\begin{proof}
For any $X\in\Gamma\left(\cal H\right)$ we have, using
\eqref{sasaki},
\begin{align*}
\left({\cal L}_{\xi_2}\phi_1\right)X &={\nabla}_{\xi_2}\left(\phi_1 X\right) -
{\nabla}_{\phi_1 X}\xi_2 - \phi_1{\nabla}_{\xi_2}X + \phi_1{\nabla}_{X}\xi_2\\
    &= \left({\nabla}_{\xi_2}\phi_1\right)X+\phi_2\phi_1 X-\phi_1\phi_2 X= -2\phi_3 X.
\end{align*}
Moreover, we have
\begin{gather*}
\left({\cal
L}_{\xi_2}\phi_1\right)\xi_1=\left[\xi_2,\phi_1\xi_1\right]-\phi_1\left[\xi_2,\xi_1\right]=2\phi_1\xi_3=-2\xi_3=-2\phi_3\xi_1,\\
\left({\cal
L}_{\xi_2}\phi_1\right)\xi_2=\left[\xi_2,\phi_1\xi_2\right]-\phi_1\left[\xi_2,\xi_2\right]=\left[\xi_2,\xi_3\right]=2\xi_1=-2\phi_3\xi_2,\\
\left({\cal
L}_{\xi_2}\phi_1\right)\xi_3=\left[\xi_2,\phi_1\xi_3\right]-\phi_1\left[\xi_2,\xi_3\right]=-\left[\xi_2,\xi_2\right]-2\phi_1\xi_1=0=-2\phi_3\xi_3,
\end{gather*}
from which we conclude ${\cal L}_{\xi_2}\phi_1=-2\phi_3$.
Similarly one can prove ${\cal L}_{\xi_3}\phi_1=2\phi_2$. Finally,
${\cal L}_{\xi_1}\phi_1=0$ holds because
$\left(\phi_1,\xi_1,\eta_1,g\right)$ is a Sasakian structure. The
other equalities in \eqref{derivatalie1} can be proved in an analogous way.
 We now prove \eqref{derivatalie2}. For any horizontal vector field $X$ we have
\begin{equation*}
\left({\cal L}_{\xi_\alpha}\phi_\beta\right)X = \nabla_{\xi_\alpha}\left(\phi_\beta X\right) - \nabla_{\phi_\beta X}\xi_\alpha - \phi_\beta\left(\nabla_{\xi_\alpha}X-\nabla_{X}\xi_\alpha\right) = \left(\nabla_{\xi_\alpha}\phi_\beta\right)X = 0
\end{equation*}
and, by using \eqref{3-sasaki} and \eqref{commutatore}, $\left({\cal
L}_{\xi_\alpha}\phi_\beta\right)\xi_\gamma=\left[\xi_\alpha,\phi_\beta
\xi_\gamma\right]-\phi_\beta\left[\xi_\alpha,\xi_\gamma\right] = 0$.
\end{proof}
\\
\\
A common property of $3$-Sasakian and $3$-cosymplectic manifolds
is stated in the following lemma.

\begin{lemma}\label{basici}
Let $M$ be a $3$-Sasakian or $3$-cosymplectic manifold. Then,
for any horizontal vector field $X$, $\left[X,\xi_\alpha\right]$ is still horizontal.
\end{lemma}
\begin{proof}
$\eta_\beta\left(\left[X,\xi_\alpha\right]\right)= X\left(\eta_\beta\left(\xi_\alpha\right)\right)-\xi_\alpha\left(\eta_\beta\left(X\right)\right)-2d\eta_\beta\left(X,\xi_\alpha\right)=-2d\eta_\beta\left(X,\xi_\alpha\right)$, for any $\beta\in\left\{1,2,3\right\}$.
Now, if the structure is $3$-cosymplectic $d\eta_\beta=0$ and if it is $3$-Sasakian
$d\eta_\beta\left(X,\xi_\alpha\right)=g\left(X,\phi_\beta\xi_\alpha\right)=\sum_{\gamma=1}^3\epsilon_{\alpha\beta\gamma}\eta_\gamma\left(X\right)=0$
since $X$ is horizontal.
\end{proof}

\bigskip

Now we attach a canonical connection to any manifold $M^{4n+3}$ with an almost contact
metric $3$-structure $\left(\phi_\alpha,\xi_\alpha,\eta_\alpha,g\right)$ in the following way. We set
\begin{equation}\label{canonica}
\tilde{\nabla}_{X}Y=\left({\nabla}_{X}Y\right)^h, \ \
\tilde{\nabla}_{\xi_\alpha}Y=\left[\xi_\alpha,Y\right], \ \
\tilde{\nabla}\xi_\alpha=0,
\end{equation}
for all $X,Y\in\Gamma\left(\cal H\right)$, where
$\left({\nabla}_{X}Y\right)^h$ denotes the horizontal component of
the Levi Civita connection. In the following proposition we start
the study of the properties of this connection.

\begin{proposition}\label{metric}
Let $\left(M^{4n+3},\phi_\alpha,\xi_\alpha,\eta_\alpha,g\right)$
be an almost $3$-contact metric manifold. Then the $1$-forms
$\eta_1, \eta_2, \eta_3$ are $\tilde{\nabla}$-parallel if and only
if $d\eta_\alpha\left(X,\xi_\beta\right)=0$ for any $X\in\Gamma\left(\cal H\right)$ and any
$\alpha,\beta\in\left\{1,2,3\right\}$. Furthermore $\tilde{\nabla}$ is a
metric connection with respect to $g$ if and only if each
$\xi_\alpha$ is Killing.
\end{proposition}
\begin{proof}
Since $\tilde{\nabla}_XY\in\Gamma\left(\cal H\right)$ for any
$X,Y\in\Gamma\left(\cal H\right)$, we have
$(\tilde{\nabla}_{X}\eta_{\alpha})Y=0$ for all
 $X,Y\in\Gamma\left(\cal H\right)$;
moreover, from $\tilde{\nabla}\xi_\beta=0$ and
$\eta_\alpha\left(\xi_\beta\right)=\delta_{\alpha\beta}$ it
follows also $(\tilde{\nabla}_E\eta_\alpha)\xi_\beta=0$, for all
$E\in\Gamma\left(TM\right)$ and
$\alpha,\beta\in\left\{1,2,3\right\}$. So $\eta_\alpha$ is
$\tilde{\nabla}$-parallel if and only if
$(\tilde{\nabla}_{\xi_\beta}\eta_\alpha)X=0$ for all
$\beta\in\left\{1,2,3\right\}$, i.e. if and only if
$\eta_\alpha\left(\left[\xi_\beta,X\right]\right)=0$ and this is
equivalent to requiring that
$d\eta_\alpha\left(X,\xi_\beta\right)=0$. Now we prove the second
part of the proposition. Firstly, we note that $\left({\cal L}_{\xi_\alpha}g\right)\left(\xi_\beta,\xi_\gamma\right) = -g\left(\left[\xi_\alpha,\xi_\beta\right],\xi_\gamma\right)-g\left(\xi_\beta,\left[\xi_\alpha,\xi_\gamma\right]\right)= -2\sum_{\delta=1}^3 \left(\epsilon_{\alpha\beta\delta} g\left(\xi_\delta,\xi_\gamma\right)
+ \epsilon_{\alpha\gamma\delta}g\left(\xi_\beta,\xi_\delta\right)\right) = -2\left(\epsilon_{\alpha\beta\gamma}+\epsilon_{\alpha\gamma\beta}\right)=0$,
and, by Lemma \ref{basici}, $\left({\cal
L}_{\xi_\alpha}g\right)\left(X,\xi_\beta\right)=\xi_\alpha\left(g\left(X,\xi_\beta\right)\right)-g\left(\left[\xi_\alpha,X\right],\xi_\beta\right)-g\left(X,2\epsilon_{\alpha\beta\gamma}\xi_\gamma\right)=0$
for $X\in \Gamma\left(\cal H\right)$. Next, we observe that for
all horizontal vector fields $X$, $Y$, $Z$, we have
\begin{align*}
(\tilde{\nabla}_{Z}g)\left(X,Y\right)&=Z\left(g\left(X,Y\right)\right)-g((\nabla_{Z}X)^h,Y)-g(X,(\nabla_{Z}Y)^h)\\
&=Z\left(g\left(X,Y\right)\right)-g\left(\nabla_{Z}X,Y\right)-g\left(X,\nabla_{Z}Y\right)=0.
\end{align*}
Moreover, clearly,
$(\tilde{\nabla}_{Z}g)\left(X,\xi_\alpha\right)=0$. Finally, $(\tilde{\nabla}_{E}g)\left(\xi_\alpha,\xi_\beta\right)=0$ for any $E\in\Gamma\left(TM\right)$ and any $\alpha,\beta\in\left\{1,2,3\right\}$.
So $g$ is $\tilde{\nabla}$-parallel if and only if
$(\tilde{\nabla}_{\xi_\alpha}g)\left(X,Y\right)=0$ for any
$X,Y\in\Gamma\left(\cal H\right)$ and for all
$\alpha\in\left\{1,2,3\right\}$. But, as $\tilde{\nabla}{\xi_\alpha}=0$, we have the equality
\begin{equation*}
(\tilde{\nabla}_{\xi_\alpha}g)\left(X,Y\right)=\xi_\alpha\left(g\left(X,Y\right)\right)-g\left(\left[{\xi_\alpha},X\right],Y\right)-g\left(X,\left[{\xi_\alpha},Y\right]\right)=\left({\cal{L}}_{\xi_\alpha}g\right)\left(X,Y\right)
\end{equation*}
from which we get the assertion.
\end{proof}
\\

In general the canonical connection $\tilde{\nabla}$ is not torsion free. Indeed we have the following result.

\begin{proposition}\label{torsion}
Let $\left(M,\phi_\alpha,\xi_\alpha,\eta_\alpha,g\right)$
be an almost $3$-contact metric manifold. Then the torsion tensor
field $\tilde{T}$ of $\tilde{\nabla}$ is given by
\begin{equation*}
 \tilde{T}\left(X,Y\right)=2\sum_{\alpha=1}^{3}
d\eta_\alpha\left(X,Y\right)\xi_\alpha, \ \
\tilde{T}\left(X,\xi_\alpha\right)=0, \ \
\tilde{T}\left(\xi_\alpha,\xi_\beta\right)=\left[\xi_\beta,\xi_\alpha\right],
\end{equation*}
for all $X, Y\in\Gamma\left(\cal H\right)$ and for all
$\alpha\in\left\{1,2,3\right\}$.
\end{proposition}
\begin{proof}
For any horizontal vector fields $X$, $Y$ we have
\begin{align*}
\tilde{T}\left(X,Y\right)&=({\nabla}_{X}Y-{\nabla}_{Y}X-\left[X,Y\right])^h-\left[X,Y\right]^v\\
&=\left(T\left(X,Y\right)\right)^h-\sum_{\alpha=1}^{3}g\left(\left[X,Y\right],\xi_\alpha\right)\xi_\alpha\\
&=2\sum_{\alpha=1}^{3}d\eta_{\alpha}\left(X,Y\right)\xi_\alpha.
\end{align*}
Moreover, it follows from $\left(\ref{canonica}\right)$ that  $\tilde{T}\left(\xi_\alpha,X\right)=\left[\xi_\alpha,X\right]-\left[\xi_\alpha,X\right]=0$.
Finally, for all $\alpha,\beta\in\left\{1,2,3\right\}$, we have easily
$\tilde{T}\left(\xi_\alpha,\xi_\beta\right)=-\left[\xi_\alpha,\xi_\beta\right]=\left[\xi_\beta,\xi_\alpha\right]$.
\end{proof}

\begin{corollary}
Let $\left(M,\phi_\alpha,\xi_\alpha,\eta_\alpha,g\right)$ be an almost $3$-contact metric manifold such that the $1$-forms
$\eta_1, \eta_2, \eta_3$ are $\tilde{\nabla}$-parallel. Then, the distribution
spanned by $\xi_1$, $\xi_2$ and $\xi_3$ is
integrable if and only if
$\tilde{T}\left(E,F\right)=2\sum_{\alpha=1}^{3}
d\eta_\alpha\left(E,F\right)\xi_\alpha$ for all
$E,F\in\Gamma\left(TM\right)$. 
\end{corollary}
\begin{proof}
From the equality
$\left[\xi_\beta,\xi_\alpha\right]^v=\sum_{\gamma=1}^{3}\eta_\gamma\left(\left[\xi_\beta,\xi_\alpha\right]\right)\xi_\gamma$
it  follows that if the distribution spanned by $\xi_1$, $\xi_2$
and $\xi_3$ is integrable, then
$\tilde{T}\left(\xi_\alpha,\xi_\beta\right)=
\sum_{\gamma=1}^{3}\eta_\gamma\left(\left[\xi_\beta,\xi_\alpha\right]\right)\xi_\gamma
= 2\sum_{\gamma=1}^{3} d\eta_\gamma
\left(\xi_\alpha,\xi_\beta\right) \xi_\gamma$. The converse is
trivial.
\end{proof}
\\

Actually, the requirement that the Reeb vector fields are parallel, together with propositions \ref{metric} and \ref{torsion} uniquely characterise the connection $\tilde{\nabla}$. This is shown in the
following theorem.

\begin{theorem}\label{connessione}
 Let  $\left(M,\phi_\alpha,\xi_\alpha,\eta_\alpha,g\right)$
be an almost $3$-contact metric manifold. Then there exists a
unique connection $\tilde{\nabla}$ on $M$ satisfying the following
properties:
\begin{description}
 \item[(i)] $\tilde{\nabla}\xi_1=\tilde{\nabla}\xi_2=\tilde{\nabla}\xi_3=0$,
 \item[(ii)] $(\tilde{\nabla}_Z g)\left(X,Y\right)=0$,\quad for all
$X,Y,Z\in\Gamma\left(\cal H\right)$,
 \item[(iii)] $ \tilde{T}\left(X,Y\right)=2\sum_{\alpha=1}^3
 d\eta_\alpha\left(X,Y\right)\xi_\alpha$ and $\tilde{T}\left(X,\xi_\alpha\right)=0$,\quad for all
 $X,Y\in\Gamma\left(\cal H\right)$.
\end{description}
Furthermore, if $M$ is $3$-Sasakian, then for all
$E,F\in\Gamma\left(TM\right)$
\begin{equation}\label{uno1}
(\tilde{\nabla}_{E}\phi_\alpha)F=-\sum_{\beta,\gamma=1}^{3}\epsilon_{\alpha\beta\gamma}\left(\eta_\beta\left(E\right)\phi_\gamma
 F^h-\eta_\gamma\left(E\right)\phi_\beta F^h\right);
\end{equation}
if $M$ is $3$-cosymplectic,  then the connection $\tilde{\nabla}$
coincides with the Levi Civita connection and in particular we
have, for each $\alpha\in\left\{1,2,3\right\}$,
$\tilde{\nabla}\phi_\alpha=0$.
\end{theorem}
\begin{proof}
The connection defined by \eqref{canonica} satisfies the properties (i)--(iii).
Thus we have only to prove the uniqueness of such a connection. Let
$\hat{\nabla}$ be any connection on $M$ verifying the
properties (i)--(iii). From (i) we get
$\hat{\nabla}\xi_\alpha=0=\tilde{\nabla}\xi_\alpha$, and, from
(iii),
$0=\hat{T}\left(\xi_\alpha,X\right)=\hat{\nabla}_{\xi_\alpha}X-\hat{\nabla}_{X}\xi_\alpha-\left[\xi_\alpha,X\right]=\hat{\nabla}_{\xi_\alpha}X-\left[\xi_\alpha,X\right]$,
which implies that
$\hat{\nabla}_{\xi_\alpha}X=\left[\xi_\alpha,X\right]=\tilde{\nabla}_{\xi_\alpha}X$
for all $X\in\Gamma\left(\cal H\right)$. Thus we have only to
verify that $\hat{\nabla}_{X}Y=\tilde{\nabla}_{X}Y$ for all
$X,Y\in\Gamma\left(\cal H\right)$, that is
$\hat{\nabla}_{X}Y=({\nabla}_{X}Y)^h$ for all
$X,Y\in\Gamma\left(\cal H\right)$. In order to check this
equality, we define another connection on $M$, by setting
\begin{equation*}
\bar{\nabla}_{E}F:=\left\{%
\begin{array}{ll}
    \hat{\nabla}_{E}F+({\nabla}_{E}F)^v, & \hbox{ for $E,F\in\Gamma\left(\cal H\right)$;} \\
    {\nabla}_{E}F, & \hbox{ for $E\in\Gamma\left({\cal H}^\perp\right)$ and
    $F\in\Gamma\left(TM\right)$;}\\
    {\nabla}_{E}F, & \hbox{ for $E\in\Gamma\left(TM\right)$ and $F\in\Gamma\left({\cal H}^\perp\right)$,}
\end{array}%
\right.
\end{equation*}
where $\left({\nabla}_EF\right)^v$ denotes the vertical component
of the Levi Civita covariant derivative. If we prove that
$\bar{\nabla}$ coincides with the Levi Civita connection, then we
will conclude that for all $X,Y\in\Gamma\left(\cal H\right)$
${\nabla}_XY=\bar{\nabla}_XY=\hat{\nabla}_XY+({\nabla}_XY)^v$,
from which $\hat{\nabla}_{X}Y=({\nabla}_XY)^h$. Firstly, note that
for all $X,Y\in\Gamma\left(\cal H\right)$, using the definition of
the Levi Civita connection $\nabla$ we have
\begin{align*}
\bar{\nabla}_{X}Y &= \hat{\nabla}_XY+\sum_{\alpha=1}^3 g\left({\nabla}_X Y,\xi_\alpha\right)\xi_\alpha\\
&= \hat{\nabla}_X Y - \frac{1}{2}\sum_{\alpha=1}^3 \left(\xi_\alpha\left(g\left(X,Y\right)\right) - g\left(\left[\xi_\alpha,X\right],Y\right)- g\left(\left[\xi_\alpha,Y\right],X\right) - g\left(\left[X,Y\right],\xi_\alpha\right)\right)\xi_\alpha\\
&= \hat{\nabla}_{X}Y + \frac{1}{2}\sum_{\alpha=1}^3\left(-\left({\cal L}_{\xi_\alpha}g\right)\left(X,Y\right) + \eta_\alpha\left(\left[X,Y\right]\right)\right)\xi_\alpha\\
&=\hat{\nabla}_XY-\sum_{\alpha=1}^3\left(\frac{1}{2}\left({\cal L}_{\xi_\alpha}g\right)\left(X,Y\right) + d\eta_\alpha\left(X,Y\right)\right)\xi_\alpha.
\end{align*}
Now we prove that the connection $\bar{\nabla}$ is metric and
torsion free. For all $X,Y,Y'\in\Gamma\left(\cal H\right)$
\begin{equation*}
(\bar{\nabla}_{X}g)\left(Y,Y'\right)=X\left(g\left(Y,Y'\right)\right)-g(\hat{\nabla}_XY,Y')-g(Y,\hat{\nabla}_XY')=(\hat{\nabla}_{X}g)\left(Y,Y'\right)=0
\end{equation*}
by the preceding equality and the condition (ii). Next, by using (iii), we obtain
$\bar{T}\left(X,Y\right)=\hat{T}\left(X,Y\right)-2\sum_{\alpha=1}^3d\eta_\alpha\left(X,Y\right)\xi_\alpha=0$.
Thus $\bar{\nabla}$ coincides with the Levi Civita connection of
$M$ and this implies that $\hat{\nabla}=\tilde{\nabla}$.

Now we prove the second part of the theorem. Assume that
$M$ is $3$-Sasakian. Then for any $X,Y\in\Gamma\left({\cal
H}\right)$, using \eqref{sasaki} and the fact that $\nabla g = 0$ we have
\begin{align*}
(\tilde{\nabla}_{X}\phi_1)Y &= (\nabla_X \phi_1)Y - \sum_{\alpha=1}^3 g\left(\nabla_X \left(\phi_1 Y\right),\xi_\alpha\right) \xi_\alpha +  \phi_1 \left(\sum_{\alpha=1}^3 g\left(\nabla_X Y, \xi_\alpha\right) \xi_\alpha \right)\\
&=  g\left(X,Y\right)\xi_1 - \eta_1\left(Y\right)X + \sum_{\alpha=1}^3 g\left(\phi_1 Y,\nabla_X \xi_\alpha \right) \xi_\alpha + g\left(\nabla_X Y, \xi_2\right) \xi_3 - g\left(\nabla_X Y, \xi_3\right) \xi_2\\
&= - g\left(\phi_1 Y, \phi_2 X\right) \xi_2 - g\left(\phi_1 Y, \phi_3 X\right) \xi_3 + g\left( Y, \phi_2 X\right) \xi_3 - g\left(Y, \phi_3 X\right) \xi_2\\
&=  g\left( Y, \phi_1\phi_2 X\right) \xi_2 + g\left( Y, \phi_1\phi_3 X\right) \xi_3 + g\left( Y, \phi_2 X\right) \xi_3 - g\left(Y, \phi_3 X\right) \xi_2 =0.
\end{align*}
Moreover, for any $\alpha,\beta,\gamma\in\left\{1,2,3\right\}$,
$(\tilde{\nabla}_{E}\phi_\beta)\xi_\gamma=\tilde{\nabla}_{E}\left(\phi_\beta
\xi_\gamma\right) - \phi_\beta\tilde{\nabla}_{E}\xi_\gamma =
\sum_{\alpha=1}^3
\epsilon_{\alpha\beta\gamma}\tilde{\nabla}_{E}\xi_\alpha=0$.
Finally, for any $X\in\Gamma\left(\cal H\right)$
\begin{equation*}
(\tilde{\nabla}_{\xi_1}\phi_1)X=\tilde{\nabla}_{\xi_1}\left(\phi_1
X\right)-\phi_1\tilde{\nabla}_{\xi_1}X=\left[\xi_1,\phi_1
X\right]-\phi_1\left[\xi_1,X\right]=\left({\cal
L}_{\xi_1}\phi_1\right)X,
\end{equation*}
so $(\tilde{\nabla}_{\xi_1}\phi_1)X=\left({\cal
L}_{\xi_1}\phi_1\right)X$.  Similarly, one can find
$(\tilde{\nabla}_{\xi_2}\phi_1)X=\left({\cal
L}_{\xi_2}\phi_1\right)X$ and
$(\tilde{\nabla}_{\xi_3}\phi_1)X=\left({\cal
L}_{\xi_3}\phi_1\right)X$. Hence, by applying \eqref{derivatalie1},
we have $(\tilde{\nabla}_{\xi_1}\phi_1)X=0$,
$(\tilde{\nabla}_{\xi_2}\phi_1)X=-2\phi_3 X$,
$(\tilde{\nabla}_{\xi_3}\phi_1)X=2\phi_2 X$. Thus, if we decompose
any pair of vector fields $E,F\in\Gamma\left(TM\right)$ in their
horizontal and vertical parts, $E=E^h+\sum_{\alpha=1}^{3}\eta_\alpha\left(E\right)\xi_\alpha$ and
$F=F^h+\sum_{\alpha=1}^{3}\eta_\alpha\left(F\right)\xi_\alpha$, we have
\begin{align*}
(\tilde{\nabla}_E\phi_1)F&=\sum_{\alpha=1}^3(\tilde{\nabla}_{\eta_\alpha\left(E\right)\xi_\alpha}\phi_1)F^h\\
&=\eta_2\left(E\right)(\tilde{\nabla}_{\xi_2}\phi_1)F^h+\eta_3\left(E\right)(\tilde{\nabla}_{\xi_3}\phi_1)F^h\\
&=-2\eta_2\left(E\right)\phi_3 F^h + 2\eta_3\left(E\right)\phi_2 F^h.
\end{align*}
The other equations involving $\phi_2$ and $\phi_3$ can be proved
in a similar way.

Finally, let $M$ be $3$-cosymplectic. Then $\nabla_{X}Y$ is horizontal
for every $X,Y\in\Gamma\left(\cal H\right)$, since $g\left(\nabla_{X}Y,\xi_\alpha\right) = -g\left(Y,\nabla_{X}\xi_\alpha\right)=0$ for all $\alpha\in\left\{1,2,3\right\}$.
Hence, $\nabla_{X}Y=(\nabla_XY)^h=\tilde{\nabla}_XY$.  Moreover, $\nabla\xi_\alpha=0=\tilde{\nabla}\xi_\alpha$. Finally,
$\nabla_{\xi_\alpha}X = \nabla_X\xi_\alpha - \left[X,\xi_\alpha\right] = \left[\xi_\alpha,X\right] = \tilde{\nabla}_{\xi_\alpha}X$. We conclude that $\nabla=\tilde{\nabla}$.
\end{proof}
\\

In the next proposition we analyse the curvature of the canonical connection
$\tilde{\nabla}$ in a $3$-Sasakian manifold.

\begin{proposition}\label{curvatura1}
Let $\left(M^{4n+3},\phi_\alpha,\xi_\alpha,\eta_\alpha,g\right)$
be a $3$-Sasakian manifold. Then the curvature
tensor of $\tilde{\nabla}$ verifies $\tilde{R}_{EF}\xi_\alpha=0$,
 $\tilde{R}_{\xi_\alpha\xi_\beta}=0$ and $\tilde{R}_{X\xi_\alpha}=0$ for all
$E,F\in\Gamma\left(TM\right)$, $X\in\Gamma\left(\cal H\right)$ and $\alpha,\beta\in\left\{1,2,3\right\}$. Moreover, for all $X,Y,Z\in\Gamma\left(\cal H\right)$,
\begin{equation}\label{formulacurvatura}
\tilde{R}_{XY}Z=\left(R_{XY}Z\right)^h+\sum_{\alpha=1}^{3}\left(d\eta_\alpha\left(Y,Z\right)\phi_\alpha X-d\eta_\alpha\left(X,Z\right)\phi_\alpha Y\right).
\end{equation}
\end{proposition}
\begin{proof}
That $\tilde{R}_{EF}\xi_\alpha=0$ is obvious since
$\tilde{\nabla}\xi_\alpha=0$. Next, for any
$\alpha,\beta\in\left\{1,2,3\right\}$,
\begin{align*}
\tilde{R}_{\xi_\alpha\xi_\beta}E&=\tilde{\nabla}_{\xi_\alpha}\left[\xi_\beta,E\right]-\tilde{\nabla}_{\xi_\beta}\left[\xi_\alpha,E\right]-2\sum_{\gamma=1}^{3}\epsilon_{\alpha\beta\gamma}\tilde{\nabla}_{\xi_\gamma}E\\
&=\left[\xi_\alpha,\left[\xi_\beta,E\right]\right]-\left[\xi_\beta,\left[\xi_\alpha,E\right]\right]-2\sum_{\gamma=1}^{3}\epsilon_{\alpha\beta\gamma}\left[\xi_\gamma,E\right]\\
&=\left[\left[\xi_\alpha,E\right],\xi_\beta\right]+\left[\left[E,\xi_\beta\right],\xi_\alpha\right]+\left[\left[\xi_\beta,\xi_\alpha\right],E\right]=0
\end{align*}
by the Jacobi identity. Moreover, since the distribution
$\left\langle\xi_1,\xi_2,\xi_3\right\rangle$ is integrable and
each $\xi_\alpha$ is Killing,  this distribution defines a
Riemannian foliation of $M^{4n+3}$, which can be described, at
least locally, by a family of Riemannian submersions. Note that
$\tilde{\nabla}$ can be interpreted as the lift of the
Levi Civita connection of the space of leaves. If $X, Y$ are
(local) basic vector fields with respect to such a given submersion, then
\begin{equation*}
\tilde{R}_{X\xi_\alpha}Y={\tilde{\nabla}_{X}}\tilde{\nabla}_{\xi_\alpha}Y-{\tilde{\nabla}_{\xi_\alpha}}\tilde{\nabla}_{X}Y-\tilde{\nabla}_{\left[X,\xi_\alpha\right]}Y=\tilde{\nabla}_{X}\left[\xi_\alpha,Y\right]-[\xi_\alpha,\tilde{\nabla}_XY]-\tilde{\nabla}_{\left[X,\xi_\alpha\right]}Y=0,
\end{equation*}
since
$\left[\xi_\alpha,Y\right]=[\xi_\alpha,\tilde{\nabla}_XY]=\left[X,\xi_\alpha\right]=0$
because, as $X,Y$ and  $\tilde{\nabla}_{X}Y$ are basic, these
brackets are vertical and, by Lemma \ref{basici}, also horizontal,
hence they vanish. It remains to prove \eqref{formulacurvatura}.
We have
\begin{align*}
\tilde R_{XY}Z&=\left(\nabla_X\tilde\nabla_Y Z\right)^h-\left(\nabla_Y\tilde\nabla_X Z\right)^h-\tilde\nabla_{\left[X,Y\right]^h}Z-\tilde\nabla_{\sum_{\alpha=1}^3\eta_\alpha\left(\left[X,Y\right]\right)\xi_\alpha}Z\\
&=\left(\nabla_X\left(\nabla_YZ-\sum_{\alpha=1}^3\eta_\alpha\left(\nabla_YZ\right)\xi_\alpha\right)\right)^h-\left(\nabla_Y\left(\nabla_XZ-\sum_{\alpha=1}^3\eta_\alpha\left(\nabla_XZ\right)\xi_\alpha\right)\right)^h\\
&\quad-\left(\nabla_{\left[X,Y\right]^h}Z\right)^h-\sum_{\alpha=1}^3\eta_\alpha\left(\left[X,Y\right]\right)\left[\xi_\alpha,Z\right]\\
&=\left(R_{XY}Z\right)^h+\sum_{\alpha=1}^3\left(\eta_\alpha\left(\nabla_YZ\right)\phi_\alpha X-\eta_\alpha\left(\nabla_XZ\right)\phi_\alpha Y\right)
\end{align*}
from which \eqref{formulacurvatura} follows.
\end{proof}
\\

We will now show that the Ricci curvature of every $3$-cosymplectic manifold
vanishes. This result is a consequence of the projectability of $3$-cosymplectic manifolds
onto hyper-K\"{a}hlerian manifolds which is stated in the following theorem.

\begin{theorem}\label{hyperkahler}
 Every regular $3$-cosymplectic structure projects onto a
hyper-K\"{a}hlerian  structure.
\end{theorem}
\begin{proof}
Since the foliation ${\cal F}_3$ is regular, it is defined by a
global submersion $f$ from $M^{4n+3}$ to the space of leaves
$M'^{4n}=M^{4n+3}/{{\cal F}_3}$. Then the Riemannian metric $g$
projects to a Riemannian metric $G$ on $M'^{4n}$ because each
$\xi_\alpha$ is Killing. Moreover, by \eqref{derivatalie2}, the
tensor fields $\phi_1$, $\phi_2$, $\phi_3$ project to three tensor
fields $J_1$, $J_2$, $J_3$ on $M'^{4n}$ and it is easy to check
that $J_\alpha
J_\beta=\sum_{\gamma=1}^{3}\epsilon_{\alpha\beta\gamma}J_\gamma-\delta_{\alpha\beta}I$.
In fact $\left(J_\alpha, G\right)$ are Hermitian structures which
are integrable because $N_\alpha=0$.
\end{proof}

\begin{remark}\label{locale}
\emph{Without the assumption of the regularity, Theorem
\ref{hyperkahler} still holds, but \emph{locally}, in the sense
that there exists a family of submersions $f_i$ from open subsets
$U_i$ of $M^{4n+3}$ to a $4n$-dimensional manifold $M'^{4n}$, with
$\left\{U_i\right\}_{i\in I}$ an open covering of $M^{4n+3}$, such
that the $3$-cosymplectic structure
$\left(\phi_\alpha,\xi_\alpha,\eta_\alpha,g\right)$ projects under
$f_i$ to a hyper-K\"{a}hlerian structure on $M'^{4n}$.}
\end{remark}

\begin{corollary}
Every $3$-cosymplectic manifold is Ricci-flat.
\end{corollary}
\begin{proof}
According to Remark \ref{locale}, let $f_i$ be a local submersion
from the $3$-cosymplectic manifold $M^{4n+3}$ to the
hyper-K\"{a}hlerian manifold $M'^{4n}$. Since $f_i$ is a
Riemannian submersion, we can apply a well-known formula which
relates the Ricci tensors and  of $M^{4n+3}$ and $M'^{4n}$ (cf.
\cite{pastore}): for any $X, Y$ basic vector fields
\begin{align}\label{ricci}
    \textrm{Ric}\left(X,Y\right)&=\textrm{Ric}'\left(f_{i_\ast}X,f_{i_\ast}Y\right) + \frac{1}{2}\left(g\left(\nabla_{X}N,Y\right)+g\left(\nabla_{Y}N,X\right)\right) \nonumber\\
&\quad-2\sum_{i=1}^{n}g\left(A_{X}X_i,A_{Y}X_i\right)-\sum_{\alpha=1}^{3}g\left(T_{\xi_\alpha}X,T_{\xi_\beta}Y\right),
\end{align}
where $\left\{X_1,\ldots,X_{4n},\xi_1,\xi_2,\xi_3\right\}$ is a
local orthonormal basis with each $X_i$ basic, $A$ and $T$ are the
O'Neill tensors associated to $f_i$, and $N$ is the local vector
field on $M^{4n+3}$ given by
$N=\sum_{\alpha=1}^{3}T_{\xi_\alpha}\xi_\alpha$. Note that, since
the horizontal distribution is integrable, $A\equiv 0$, and by
$\nabla\xi_\alpha=0$ we get
$T_{\xi_\alpha}\xi_\alpha=\left(\nabla_{\xi_\alpha}\xi_\alpha\right)^h=0$,
$T_{\xi_\alpha}Z=\left(\nabla_{\xi_\alpha}Z\right)^v=\left(\nabla_{Z}\xi_\alpha+\left[\xi_\alpha,Z\right]\right)^v=0$.
Hence the formula \eqref{ricci} reduces to
\begin{equation*}
\textrm{Ric}\left(X,Y\right)=\textrm{Ric}'\left(f_{i_\ast}X,f_{i_\ast}Y\right).
\end{equation*}
But $\textrm{Ric}'\left(X',Y'\right)=0$ for all
$X',Y'\in\Gamma\left(TM\right)$, because $M'^{4n}$ is
hyper-K\"{a}hlerian. Hence $\textrm{Ric}=0$ in the horizontal
subbundle $\cal H$. Finally, it is easy to check that
$\textrm{Ric}\left(\xi_\alpha,\xi_\beta\right)=0$ and
$\textrm{Ric}\left(X,\xi_\beta\right)=0$ for any
$X\in\Gamma\left(\cal H\right)$.
\end{proof}

\section{The Darboux Theorem}
Let $M^{4n+3}$ be a manifold endowed with an almost contact metric
$3$-structure $\left(\phi_\alpha,\xi_\alpha,\eta_\alpha,g\right)$.
We denote by $\Phi_\alpha^\flat\colon X \mapsto
\Phi_\alpha(X,\cdot\,)$ the musical isomorphisms induced by the
fundamental $2$-forms $\Phi_\alpha$ between horizontal vector
fields and vertical $1$-forms. Their inverses will be denoted by
$\Phi_\alpha^\sharp$. We  also denote by $g^{\flat}_\mathcal{H}$
the musical isomorphism induced by the metric between horizontal
vector fields and vertical $1$-forms, and by
$\phi^{\mathcal{H}}_\alpha\colon \mathcal{H} \rightarrow
\mathcal{H}$ the isomorphisms induced  by the endomorphisms
$\phi_\alpha\colon TM \rightarrow TM$.
\begin{lemma}\label{antonio}
In any almost $3$-contact metric manifold, the following formulas hold, for each $\alpha\in\left\{1,2,3\right\}$,
\begin{equation}\label{formulaantonio}
g^{\flat}_\mathcal{H} = \Phi_\alpha^\flat\circ\phi^{\mathcal{H}}_\alpha, \qquad
\phi^{\mathcal{H}}_\alpha = -\frac{1}{2} \sum_{\beta,\gamma=1}^3\epsilon_{\alpha\beta\gamma}\Phi_{\beta}^{\sharp}\circ
\Phi_{\gamma}^\flat.
\end{equation}
\end{lemma}
\begin{proof}
From $\Phi_\alpha\left(X,Y\right)=g\left(X,\phi_\alpha Y\right)$
we have
$-\Phi_{\alpha}^{\flat}=g^{\flat}_\mathcal{H}\circ\phi^{\mathcal{H}}_{\alpha}$.
It follows that
\begin{equation}\label{equazione1}
g^{\flat}_\mathcal{H}=\Phi_{\alpha}^{\flat}\circ\phi^{\mathcal{H}}_{\alpha},
\end{equation}
since $\phi_{\alpha}^{2}X=-X+\eta_{\alpha}(X)\xi_\alpha=-X$ for every $X\in\Gamma\left(\mathcal{H}\right)$.
We now prove the second formula of \eqref{formulaantonio}.
Since the equation \eqref{equazione1} holds for each $\alpha\in\left\{1,2,3\right\}$, we get
\begin{equation}\label{alfabeta}
\phi^{\mathcal{H}}_{\beta}\circ\phi^{\mathcal{H}}_{\gamma}  = -\Phi_{\beta}^{\sharp} \circ \Phi_{\gamma}^{\flat},
\end{equation}
for each $\beta,\gamma\in\left\{1,2,3\right\}$. Moreover, in view of  \eqref{3-sasaki}, we have $\phi^{\mathcal{H}}_{\beta}\circ\phi^{\mathcal{H}}_{\gamma} = \sum_{\alpha=1}^3\epsilon_{\alpha\beta\gamma}\phi^{\mathcal{H}}_\alpha$.
Thus we obtain $\sum_{\alpha=1}^3\epsilon_{\alpha\beta\gamma}\phi^{\mathcal{H}}_\alpha = -\Phi_{\beta}^{\sharp} \circ \Phi_{\gamma}^{\flat}$, that is
$2\phi^{\mathcal{H}}_\alpha = -\sum_{\beta,\gamma=1}^3\epsilon_{\alpha\beta\gamma}\Phi_{\beta}^{\sharp} \circ \Phi_{\gamma}^{\flat}$.
\end{proof}

\begin{corollary}\label{antonio2}
In any almost $3$-contact metric manifold, the following formula holds in
the horizontal subbundle $\mathcal{H}$,%
\[
g^{\flat}_\mathcal{H} = -\Phi_{1}^{\flat}\circ \Phi_{2}^{\sharp}\circ \Phi_{3}^{\flat}.
\]
\end{corollary}
\begin{proof}
From the two equalities in \eqref{formulaantonio} we obtain
\[
g^{\flat}_\mathcal{H} = -\frac{1}{2}\Phi_1^\flat \circ \left(\Phi_2^{\sharp}\circ\Phi_{3}^\flat - \Phi_3^{\sharp}\circ\Phi_{2}^\flat\right).
\]
On the other hand, from \eqref{alfabeta} and \eqref{3-sasaki} we obtain $\Phi_2^\flat\circ\phi^{\mathcal{H}}_3 = -\Phi_3^\flat\circ\phi^{\mathcal{H}}_2$. The claim follows.
\end{proof}
\\

Now we prove that a $3$-Sasakian manifold cannot admit any
Darboux-like coordinate system. Here for ``Darboux-like coordinate
system'' we mean local coordinates
$\left\{x_1,\ldots,x_{4n},z_1,z_2,z_3\right\}$ with respect to
which, for each $\alpha\in\left\{1,2,3\right\}$, the fundamental
$2$-forms $\Phi_\alpha=d\eta_\alpha$ have constant components and
$\xi_\alpha=a^1_\alpha \frac{\partial}{\partial z_1} +
a^2_\alpha\frac{\partial}{\partial z_2} + a^3_\alpha
\frac{\partial}{\partial z_3}$,  $a_\alpha^\beta$ being functions
depending only on the coordinates $z_1,z_2,z_3$. This is a natural
generalisation of the standard Darboux coordinates for contact
manifolds.

\begin{theorem}\label{sasakiano}
Let $\left(M^{4n+3},\phi_\alpha,\xi_\alpha,\eta_\alpha,g\right)$
be a $3$-Sasakian manifold. Then $M^{4n+3}$ does not admit any
Darboux-like coordinate system.
\end{theorem}
\begin{proof}
Let $p$ be a point of $M^{4n+3}$. Then in view of Theorem
\ref{proiezione} there exist an open neighbourhood $U$ of $p$ and
a (local) Riemannian submersion $f$ with connected fibres from $U$
onto a quaternionic K\"{a}hlerian manifold $M'^{4n}$, such that
$\ker\left(f_\ast\right)=\left\langle\xi_1,\xi_2,\xi_3\right\rangle$.
Note that the horizontal vectors with respect to $f$ are just the
vectors belonging to $\cal H$, i.e. those orthogonal to
$\xi_1,\xi_2,\xi_3$. Now, suppose by contradiction that about the point
$p$ there exists a Darboux coordinate system, that is an open
neighbourhood $V$ with local coordinates
$\left\{x_1,\ldots,x_{4n},z_1,z_2,z_3\right\}$ as above. We can
assume that $U=V$. We decompose each vector field
$\frac{\partial}{\partial x_i}$ in its horizontal and vertical
components, $\frac{\partial}{\partial x_i} = X_i +
\sum_{\alpha=1}^{3} \eta_\alpha\left(\frac{\partial}{\partial
x_i}\right)\xi_\alpha$. Note that
\begin{align}\label{dipendenzadaz}
\eta_\alpha\left(\frac{\partial}{\partial x_i}\right) &=  \frac{1}{2}\sum_{\beta,\gamma=1}^{3} \epsilon_{\alpha\beta\gamma} g\left(\frac{\partial}{\partial x_i}, \phi_\beta\xi_\gamma\right)\nonumber = \frac{1}{2}\sum_{\beta,\gamma=1}^{3} \epsilon_{\alpha\beta\gamma} d\eta_\beta\left(\frac{\partial}{\partial x_i},\xi_\gamma\right)\\
&= \frac{1}{2}\sum_{\beta,\gamma,\delta=1}^{3} \epsilon_{\alpha\beta\gamma} a^\delta_\gamma d\eta_\beta\left(\frac{\partial}{\partial x_i}, \frac{\partial}{\partial z_\delta}\right),
\end{align}
so that $\eta_\alpha\left(\frac{\partial}{\partial x_i}\right)$
are functions which do not depend on the coordinates $x_i$.
Consequently, the only eventually non-constant components of each
horizontal vector field $X_i=\frac{\partial}{\partial
x_i}-\sum_{\alpha=1}^{3} \eta_\alpha\left(\frac{\partial}{\partial
x_i}\right)\xi_\alpha$ in the holonomic basis
$\left(\frac{\partial}{\partial
x_1},\ldots,\frac{\partial}{\partial
x_{4n}},\frac{\partial}{\partial z_1},\frac{\partial}{\partial
z_2},\frac{\partial}{\partial z_3}\right)$ depend at most on the
coordinates $z_1,z_2,z_3$. Actually, for each
$i\in\left\{1,\ldots,4n\right\}$, $X_i$ is a basic vector field
with respect to the submersion $f$, thus its components do not
depend even on the fibre coordinates $z_\alpha$, hence they are
constant. For proving this it is sufficient to show that, for each
$\alpha\in\left\{1,2,3\right\}$, $\left[X_i,\xi_\alpha\right]$ is
vertical. Indeed,
\begin{equation*}
\left[X_i,\xi_\alpha\right]= \sum_{\beta=1}^{3} \frac{\partial
a^\beta_\alpha}{\partial x_i} \frac{\partial}{\partial z_\beta}
+\sum_{\beta=1}^3\left[\eta_\beta\left(\frac{\partial}{\partial
x_i}\right)\xi_\beta,\xi_\alpha\right]=\sum_{\beta=1}^3\left[\eta_\beta\left(\frac{\partial}{\partial
x_i}\right)\xi_\beta,\xi_\alpha\right]
\end{equation*}
because the functions $a^\beta_\alpha$ do not depend on the
coordinates $x_i$. Then by Corollary \ref{antonio2}
\begin{equation}
g\left(X_i,X_j\right)=-\left(d\eta_{1}^{\flat}\circ d\eta_{2}^{\sharp}\circ d\eta_{3}^{\flat}\right)\left(X_i\right)\left(X_j\right)
\end{equation}
and so the functions $g\left(X_i,X_j\right)$ are constant since
each $X_i$ has constant components and the $2$-forms
$d\eta_\alpha$ are assumed to have constant components, too. The
next step is to note that, for all
$i,j\in\left\{1,\ldots,4n\right\}$, the brackets
$\left[X_i,X_j\right]$ are vertical vector fields. We have, by
\eqref{dipendenzadaz},
\begin{align*}
\left[X_i,X_j\right]&=\left[\frac{\partial}{\partial x_i},\frac{\partial}{\partial
x_j}\right]+\sum_{\alpha,\beta=1}^{3}\left[\eta_\alpha\left(\frac{\partial}{\partial
x_i}\right)\xi_\alpha,\eta_\beta\left(\frac{\partial}{\partial
x_j}\right)\xi_\beta\right]\\
&\quad-\sum_{\alpha=1}^{3}\left[\eta_\alpha\left(\frac{\partial}{\partial
x_i}\right)\xi_\alpha,\frac{\partial}{\partial
x_j}\right]-\sum_{\beta=1}^{3}\left[\frac{\partial}{\partial
x_i},\eta_\beta\left(\frac{\partial}{\partial
x_j}\right)\xi_\beta\right]\\
&= 2\sum_{\alpha,\beta,\gamma=1}^{3} \eta_\alpha\left(\frac{\partial}{\partial
x_i}\right)\eta_\beta\left(\frac{\partial}{\partial
x_j}\right) \epsilon_{\alpha\beta\gamma} \xi_\gamma.
\end{align*}
Then, for all $i,j,k\in\left\{1,\ldots,4n\right\}$, using \eqref{canonica} and the Koszul formula for the Levi Civita covariant derivative we obtain
\begin{align*}
2g(\tilde{\nabla}_{X_i}X_j,X_k)=2g\left(\nabla_{X_i}X_j,X_k\right)&=X_i\left(g\left(X_j,X_k\right)\right)+X_j\left(g\left(X_k,X_i\right)\right)-X_k\left(g\left(X_i,X_j\right)\right)\\
&\quad-g\left(\left[X_j,X_k\right],X_i\right)+g\left(\left[X_k,X_i\right],X_j\right)+g\left(\left[X_i,X_j\right],X_k\right)=0,
\end{align*}
so that $\tilde{\nabla}_{X_i}X_j=0$. But $\tilde{\nabla}$ projects
locally to the Levi Civita connection $\nabla'$ of the
quaternionic K\"{a}hlerian manifold $M'^{4n}$ under the Riemannian
submersion $f$ so that in particular we would have that $\nabla'$
is flat and this cannot happen because the scalar curvature of
$M'^{4n}$, by Theorem \ref{proiezione}, must be strictly positive.
\end{proof}
\\

Now we prove a Darboux theorem for $3$-cosymplectic manifolds.

\begin{theorem}\label{cosimplettico}
Around each point of a flat $3$-cosymplectic manifold
$M^{4n+3}$ there are local coordinates
$\left\{x_1,\ldots,x_n,y_1,\ldots,y_n,u_1,\ldots,u_n,v_1,\ldots,v_n,z_1,z_2,z_3\right\}$
such that, for each $\alpha\in\left\{1,2,3\right\}$,
$\eta_\alpha=dz_\alpha$, $\xi_\alpha=\frac{\partial}{\partial
z_\alpha}$ and, moreover,
\begin{gather}
\Phi_1=2\sum_{i=1}^{n}\left(dx_i\wedge dy_i+du_i\wedge dv_i\right)-2dz_2\wedge dz_3,\label{Phiprima}\\
\Phi_2=2\sum_{i=1}^{n}\left(dx_i\wedge du_i-dy_i\wedge dv_i\right)+2dz_1\wedge dz_3,\\
\Phi_3=2\sum_{i=1}^{n}\left(dx_i\wedge dv_i+dy_i\wedge
du_i\right)-2dz_1\wedge dz_2,\label{Phiterza}
\end{gather}
$\phi_1$, $\phi_2$ and $\phi_3$ are represented, respectively, by the
$(4n+3)\times(4n+3)$-matrices {\small
\begin{equation}\label{matrice1}
\phi_1=
\left(%
\begin{array}{ccccccc}
  0 & -I_n & 0 & 0 & 0 & 0 & 0 \\
  I_n & 0 & 0 & 0 & 0 & 0 & 0 \\
  0 & 0 & 0 & -I_n & 0 & 0 & 0 \\
  0 & 0 & I_n & 0 & 0 & 0 & 0 \\
  0 & 0 & 0 & 0 & 0 & 0 & 0 \\
  0 & 0 & 0 & 0 & 0 & 0 & -1 \\
  0 & 0 & 0 & 0 & 0 & 1 & 0 \\
\end{array}%
\right),
\end{equation}
\begin{gather}
\phi_2=
\left(%
\begin{array}{ccccccc}\label{matrice2}
  0 & 0 & -I_n & 0 & 0 & 0 & 0 \\
  0 & 0 & 0 & I_n & 0 & 0 & 0 \\
  I_n & 0 & 0 & 0 & 0 & 0 & 0 \\
  0 & -I_n & 0 & 0 & 0 & 0 & 0 \\
  0 & 0 & 0 & 0 & 0 & 0 & 1 \\
  0 & 0 & 0 & 0 & 0 & 0 & 0 \\
  0 & 0 & 0 & 0 & -1 & 0 & 0 \\
\end{array}%
\right),\\
\label{matrice3}\phi_3=
\left(%
\begin{array}{ccccccc}
  0 & 0 & 0 & -I_n & 0 & 0 & 0 \\
  0 & 0 & -I_n & 0 & 0 & 0 & 0 \\
  0 & I_n & 0 & 0 & 0 & 0 & 0 \\
  I_n & 0 & 0 & 0 & 0 & 0 & 0 \\
  0 & 0 & 0 & 0 & 0 & -1 & 0 \\
  0 & 0 & 0 & 0 & 1 & 0 & 0 \\
  0 & 0 & 0 & 0 & 0 & 0 & 0 \\
\end{array}%
\right).
\end{gather}}
\end{theorem}
\begin{proof}
Let $p$ be a point of $M^{4n+3}$. Since $M^{4n+3}$ is flat there
exists a neighbourhood $U$ of $p$ where the curvature tensor field
vanishes identically. Moreover, one can prove by some linear
algebra that there exist horizontal vectors $e_1,\ldots,e_n$ such
that
$\{e_1,\ldots,e_n,\phi_1e_1,\ldots,\phi_1e_n,\phi_2e_1,\ldots,\phi_2e_n,$
$\phi_3e_1,\ldots,\phi_3e_n,\xi_{1_p},\xi_{2_p},\xi_{3_p}\}$ is an
orthonormal basis of $T_p M$ satisfying the equalities
\begin{gather*}
\Phi_1\left(e_i,\phi_1 e_j\right)=\delta_{ij},
\Phi_1\left(\phi_2 e_i,\phi_3 e_j\right)=\delta_{ij}, \Phi_1\left(\xi_{2_p},\xi_{3_p}\right)=-1,\\
\Phi_2\left(e_i,\phi_2 e_j\right)=\delta_{ij},
\Phi_2\left(\phi_1 e_i,\phi_3 e_i\right)=-\delta_{ij}, \Phi_2\left(\xi_{1_p},\xi_{3_p}\right)=1,\\
\Phi_3\left(e_i,\phi_3 e_j\right)=\delta_{ij}, \Phi_3\left(\phi_1
e_i,\phi_2 e_j\right)=\delta_{ij},
\Phi_3\left(\xi_{1_p},\xi_{2_p}\right)=-1,
\end{gather*}
and such that the values of the 2-forms $\Phi_\alpha$ on all the other pairs of basis vectors vanish.
Now we define $4n$ vector fields $X_i$, $Y_i$, $U_i$, $V_i$ on $U$
by  parallel transport  of the vectors $e_i$,
$\phi_1 e_i$, $\phi_2 e_i$, $\phi_3 e_i$,
$i\in\left\{1,\ldots,n\right\}$. Note that the definition is
well-posed because the parallel transport does not depend on the
curve. Since the Levi Civita connection is a metric connection and
since $\nabla\xi_\alpha=0$ we have that
$\{X_1,\ldots,X_n,Y_1,\ldots,Y_n,U_1,\ldots,U_n,$
$V_1,\ldots,V_n,\xi_1,\xi_2,\xi_3\}$ is an orthonormal frame on
$U$. Moreover by $\nabla\phi_\alpha=0$ we get that
\begin{equation}\label{fi}
Y_i=\phi_1 X_i, \ U_i=\phi_2 X_i, \ V_i=\phi_3 X_i,
\end{equation}
and by $\nabla \Phi_\alpha=0$ we have
\begin{gather}\label{prima3}
\Phi_1\left(X_i,Y_j\right)=\delta_{ij},
\Phi_1\left(U_i,V_j\right)=\delta_{ij}, \Phi_1\left(\xi_{2},\xi_{3}\right)=-1,\\
\Phi_2\left(X_i,U_j\right)=\delta_{ij},
\Phi_2\left(Y_i,V_j\right)=-\delta_{ij}, \Phi_2\left(\xi_{1},\xi_{3}\right)=1,\\
\label{ultima3} \Phi_3\left(X_i,V_j\right)=\delta_{ij},
\Phi_3\left(Y_i,U_j\right)=\delta_{ij},
\Phi_3\left(\xi_{1},\xi_{2}\right)=-1,
\end{gather}
and the values of the 2-forms $\Phi_\alpha$ on all the other pairs of vector fields
belonging to the orthonormal frame vanish.
Since the vector fields $X_i$, $Y_i$, $U_i$, $V_i$ are, by
construction, $\nabla$-parallel we have that the bracket of each pair
of these vector fields vanishes identically. This, together with
\eqref{commutatore} and the vanishing of the brackets $\left[X_i,\xi_\alpha\right]$,  $\left[Y_i,\xi_\alpha\right]$,
$\left[U_i,\xi_\alpha\right]$ and  $\left[V_i,\xi_\alpha\right]$ implies the existence of local coordinates
$\left\{x_1,\ldots,x_n,y_1,\ldots,y_n,\right.$ \allowbreak $\left. u_1,\ldots,u_n,v_1,\ldots,v_n,z_1,z_2,z_3 \right\}$
with respect to which
\begin{equation*}
X_i=\frac{\partial}{\partial x_i}, \
Y_i=\frac{\partial}{\partial y_i}, \
U_i=\frac{\partial}{\partial u_i}, \
V_i=\frac{\partial}{\partial v_i}, \
\xi_1=\frac{\partial}{\partial z_1}, \
\xi_2=\frac{\partial}{\partial z_2}, \
\xi_3=\frac{\partial}{\partial z_3}.
\end{equation*}
Now, as  the 1-forms $\eta_\alpha$ are closed, they are locally exact, and  we have (eventually reducing $U$)
$\eta_\alpha=df_\alpha$ for some functions $f_\alpha\in C^\infty\left(U\right)$,
and from the relations $\eta_\alpha\left(X_i\right)=\eta_\alpha\left(Y_i\right) = \eta_\alpha\left(U_i\right)=\eta_\alpha\left(V_i\right) = 0$,
$\eta_\alpha\left(\xi_\beta\right)=\delta_{\alpha\beta}$ it follows that
$\frac{\partial f_\alpha}{\partial x_i} = \frac{\partial f_\alpha}{\partial y_i}
= \frac{\partial f_\alpha}{\partial u_i} = \frac{\partial f_\alpha}{\partial v_i} = 0$,
$\frac{\partial f_\alpha}{\partial z_\beta} = \delta_{\alpha\beta}$.
Hence, for each $\alpha\in\left\{1,2,3\right\}$, $\eta_\alpha=dz_\alpha$.
Next, by \eqref{prima3}--\eqref{ultima3}, we get \eqref{Phiprima}--\eqref{Phiterza}.
Finally, by \eqref{fi} and by $\phi_\alpha\xi_\beta = \sum_{\gamma=1}^3
\epsilon_{\alpha\beta\gamma} \xi_\gamma$ we deduce that with respect to this coordinate
system $\phi_1$, $\phi_2$ and $\phi_3$ are represented by the matrices \eqref{matrice1},
\eqref{matrice2} and \eqref{matrice3}, respectively.
\end{proof}
\\

Arguing as in Theorem \ref{sasakiano} and taking into account that the "vertical" terms
$R_{\xi_\alpha\xi_\beta}$ and the "mixed" terms $R_{X\xi_\alpha}$ of the curvature tensor
(with $X\in\Gamma\left(\cal H\right)$) vanish,
one can prove the converse of Theorem \ref{cosimplettico}:

\begin{proposition}
Let $\left(M^{4n+3},\phi_\alpha,\xi_\alpha,\eta_\alpha,g\right)$
be a $3$-cosymplectic manifold. If each point of $M^{4n+3}$ admits
a Darboux coordinate system such that \eqref{Phiprima}--\eqref{Phiterza} of Theorem \ref{cosimplettico} hold, then
$M^{4n+3}$ is flat.
\end{proposition}

\bigskip

\begin{remark}
\emph{We conclude noting that in any almost $3$-contact metric
manifold $\left(M,\phi_\alpha,\xi_\alpha,\eta_\alpha,g\right)$
(and in particular in any hyper-contact manifold (cf.
\cite{banyaga})) the metric $g$ is uniquely determined by the
three fundamental $2$-forms $\Phi_\alpha$ and the three Reeb
vector fields $\xi_\alpha$. In particular, in the case of
$3$-Sasakian manifolds the metric is uniquely determined by the
three contact forms $\eta_\alpha$. Indeed, on the one hand, it
follows from Corollary \ref{antonio2} that}
\[g \left(X,Y\right) = -(d\eta_{1}^{\flat} \circ d\eta_{2}^{\sharp} \circ d\eta_{3}^{\flat}\left(X\right)) \left(Y\right),\]
\emph{for any $X,Y\in \Gamma\left(\cal H\right)$. On the other hand, we have $g\left(\xi_\alpha,\xi_\beta \right) = \delta_{\alpha\beta}$ and $g\left(X,\xi_\alpha\right)=\eta_\alpha\left(X\right)=0$.
This remark gives an answer to the open problem raised by Banyaga in the Remark 11 of \cite{banyaga}.}
\end{remark}

\end{document}